\documentclass[a4paper,english,12pt]{amsart}

\usepackage{amscd}
\usepackage{srcltx}
\usepackage{amsmath}
\usepackage{amsfonts}
\usepackage{verbatim}
\usepackage{amsthm}
\usepackage{amssymb, mathrsfs}
\usepackage{amscd}
\usepackage{latexsym}
\usepackage{array}
\usepackage{enumerate}
\usepackage{longtable}
\usepackage[all]{xypic}
\usepackage[latin1]{inputenc}

\ifx\pdfpagewidth\undefined 
\usepackage[T1]{fontenc}
\else 
\usepackage[OT1]{fontenc} 
\fi 
\usepackage{amsfonts,amssymb,latexsym,longtable,amsmath}

\newtheorem{thm}{Theorem}[section]
\newcommand{\bthm}{\begin{thm}} \newcommand{\ethm}{\end{thm}}
\newtheorem{prop}[thm]{Proposition}
\newcommand{\bprp}{\begin{prop}} \newcommand{\eprp}{\end{prop}}
\newtheorem{fact}[thm]{Fact}
\newcommand{\bfct}{\begin{fact}} \newcommand{\efct}{\end{fact}}
\newtheorem{prob}[thm]{Problem}
\newcommand{\bprb}{\begin{prob}} \newcommand{\eprb}{\end{prob}}
\newtheorem{question}[thm]{Question}
\newcommand{\bque}{\begin{question}} \newcommand{\eque}{\end{question}}
\newtheorem{lem}[thm]{Lemma}
\newcommand{\blem}{\begin{lem}} \newcommand{\elem}{\end{lem}}
\newtheorem{claim}[thm]{Claim}
\newcommand{\bclm}{\begin{claim}} \newcommand{\eclm}{\end{claim}}
\newtheorem{cor}[thm]{Corollary}
\newcommand{\bcor}{\begin{cor}} \newcommand{\ecor}{\end{cor}}
\newtheorem{conj}[thm]{Conjecture}
\newcommand{\bcnj}{\begin{conj}} \newcommand{\ecnj}{\end{conj}}
\theoremstyle{definition}
\newtheorem{defn}[thm]{Definition}
\newcommand{\bdfn}{\begin{defn}} \newcommand{\edfn}{\end{defn}}
\newtheorem{spec}[thm]{Specializing}
\newcommand{\bspc}{\begin{spec}} \newcommand{\espc}{\end{spec}}
\theoremstyle{remark}
\newtheorem{rem}[thm]{Remark}
\newcommand{\brem}{\begin{rem}} \newcommand{\erem}{\end{rem}}
\newtheorem{cnv}[thm]{Convention}
\newcommand{\bcnv}{\begin{cnv}} \newcommand{\ecnv}{\end{cnv}}
\newtheorem{exam}[thm]{{\bf Example }}
\newcommand{\bexm}{\begin{exam}} \newcommand{\eexm}{\end{exam}}
\newcommand{\bpf}{\begin{proof}} \newcommand{\epf}{\end{proof}}

\newcommand{\C}{\mathbb C}

\newcommand{\N}{\mathbb N}

\newcommand{\Q}{\mathbb Q}

\newcommand{\NN} {{\mathbb N}}

\renewcommand{\phi}{\varphi}
\renewcommand{\theta}{\vartheta}

\newcommand{\gep}{{\epsilon}}

\newcommand{\s}{{\sigma}}
\newcommand{\irn}{\textrm{irrep}_n(G)}

\newcommand{\irs}{\textrm{irrep}}
\newcommand{\irm}{\textrm{irrep}_m(G)}

\newcommand{\rep}{\textrm{rep}}
\newcommand{\repg}{\textrm{rep}(G)}

\newcommand{\repng}{\textrm{rep}_{n}(G)}
\newcommand{\repnh}{\textrm{rep}_{n}(H)}

\newcommand{\U}{\mathbb{U}}

\newcommand{\nbd}{neighborhood \ }

\begin{document}
\title[Bohr compactification and Chu duality of locally compact groups]{Bohr compactification and Chu duality of non-abelian locally compact groups}

\author[M. V. Ferrer]{Mar\'ia V. Ferrer}
\address{Universitat Jaume I, IMAC and Departamento de Matem\'{a}ticas,
Campus de Riu Sec, 12071 Castell\'{o}n, Spain.}
\email{mferrer@mat.uji.es}
\author[S. Hern\'andez]{S. Hern\'andez}
\address{Universitat Jaume I, IMAC and Departamento de Matem\'{a}ticas,
Campus de Riu Sec, 12071 Castell\'{o}n, Spain.}
\email{hernande@mat.uji.es}
\date{\today}

\thanks{ The first and third authors acknowledge partial support by by the Spanish Ministerio de Econom\'{i}a y Competitividad,
grant: MTM/PID2019-106529GB-I00 (AEI/FEDER, EU) and by the Universitat Jaume I, grant UJI-B2022-39}

\begin{abstract}
The Bohr compactification of an arbitrary topological group $G$ is defined as the group compactification $(bG,b)$ with the following universal property:
for every continuous homomorphism $h$ from $G$ into a compact group $K$ there is a continuous homomorphism $h^{b}$ from $bG$ into $K$
that extends $h$ in the sense that $h=h^b \circ b$. The Bohr compactification $(bG,b)$ is the unique (up to equivalence) largest compactification of $G$.
Although, for locally compact abelian groups the Bohr compactification is a big monster, for non-Abelian groups the situation is much more interesting
and one can say that all options are possible. Here we discuss the relation between the notions of duality and the Bohr compactification of
non-Abelian locally compact groups. Our main result establishes the existence of noncompact Chu reflexive groups whose partial dual spaces $\widehat G_n$,
equipped with the Fell topology, are discrete. Furthermore, we characterize when the Bohr compactification of a locally compact group is topologically isomorphic
to its Chu or unitary quasi-dual. These results fix some incorrect statements that have appeared in the literature.
\end{abstract}


\maketitle \setlength{\baselineskip}{24pt}
\setlength{\parindent}{0.5cm}

\section{Introduction}
Suppose an algebraic group $G$  is equipped with two locally compact topologies
$G_1=(G,\tau_1)$ and $G_2=(G,\tau_2)$. In case $G$ is abelian, we have that the identity mapping of ${\rm i_d}\colon G_1\to G_2$
is a topological group isomorphism if and only if there is a topological isomorphism $\Phi\colon bG_1\to bG_2$
such that $\Phi_{|G}= {\rm i_d}$, where
$bG_1$ and $bG_2$ denote the the Bohr compactifications of $G_1$ and $G_2$, respectively.
Furthermore, if $\tau_1\subsetneq \tau_2$, then $|\frac{bG_2}{bG_1}|\geq 2^\frak c$.
Therefore the Bohr compactification of a locally compact abelian group completely
characterizes its topological and algebraic structure.

This fact does not extend to non-Abelian groups, where essentially all options are possible.
This occurs because, for a non-Abelian group, the set of finite-dimensional representations no longer reflects
its overall topological and algebraic structure. In fact, there exist compact groups that are the Bohr compactification of themselves
when equipped with the discrete topology (these are known as self-Bohrifying groups).

In this project, we aim to clarify the relationship between the Bohr compactification of locally compact groups
and the Chu (or unitary) duality of these groups. For example, we characterize the conditions under which
the Bohr compactification of a locally compact group is topologically isomorphic to its Chu or unitary quasi-dual.

We also address the question of whether a group whose Chu dual is discrete must necessarily be compact.
Specifically, we establish the existence of non-compact Chu-reflexive groups whose partial dual spaces,
$\widehat{G}_n$, equipped with the Fell topology, are discrete.
Our results correct some inaccurate statements found in the literature.

\section{Preliminaries}
All topological groups are assumed to be Hausdorff.
For a (complex) Hilbert space $\mathcal{H}$ the unitary group $\mathbb U(\mathcal{H})$
of all linear isometries of $\mathcal{H}$ is equipped with the strong operator topology
(this is the topology of pointwise convergence). With this topology, $\mathbb U(\mathcal{H})$  is a topological group.
If $\mathcal{H}=\C^n$, we identify $\mathbb U(\mathcal{H})$ with the {\em unitary group $\mathbb{U}(n)$
of order $n$}, that is, the compact Lie group of all complex $n\times n$ matrices $M$
for which $M^{-1}=M^{*}$.

A {\em unitary representation} $\rho$ of the to\-po\-lo\-gi\-cal group $G$
is a continuous homomorphism $G\to \mathbb U(\mathcal{H})$, where $\mathcal{H}$ is a complex Hilbert space.
A closed linear subspace $E\subseteq \mathcal H$ is an \emph{invariant} subspace
for $\mathcal S\subseteq \mathbb U(\mathcal{H})$ if $ME\subseteq E$ for all
$M\in \mathcal S$.
If there is a closed subspace $E$ with $\{0\}\subsetneq E\subsetneq \mathcal H$
which is invariant for $\mathcal S$, then $\mathcal S$ is called
\emph{reducible}; otherwise $\mathcal S$ is \emph{irreducible}.
An \emph{irreducible representation} of $G$ is a 
unitary representation $\rho$ such that 
$\rho(G)$ is irreducible.

Two unitary representations $\rho:G\to \mathbb U(\mathcal{H}_1)$ and $\psi: G\to \mathbb U(\mathcal{H}_2)$
are {\it equivalent} ($\rho \sim \psi$)
if there exists a Hilbert space isomorphism $M:\mathcal {H}_1\to \mathcal H_2$
such that $\rho(x)=M^{-1}\psi(x)M$ for all
$x\in G$. The {\em dual object} of a topological group $G$ is the set $\widehat G$ of
equivalence classes of irreducible unitary representations of $G$. Sometimes it is more useful to consider the subset of representations with a  \emph{finite range}.
In that sense, we will denote by $\widehat{G}^f$ the subset of the equivalence classes corresponding representations with a finite range.

If $G$ is a precompact group, the Peter-Weyl Theorem (see \cite{hof_mor:compact_groups}) implies that
all irreducible unitary representation of $G$ are finite-dimensional and determine an embedding of $G$ into the product
of unitary groups $\mathbb{U}(n)$.

If $\rho:G\to \mathbb U(\mathcal{H})$ is a unitary representation, a complex-valued function $f$ on $G$ is called
a {\em function of positive type} (or {\em positive definite function}) {\em associated with} $\rho$ if there exists a vector
$v\in \mathcal{H}$ such that $f(g)=(\rho(g)v, v)$
(here $(\cdot ,\cdot)$ denotes the inner product in $\mathcal{H}$).  We denote by
$P_\rho'$ be the set of all functions of positive type associated with $\rho$.
Let $P_\rho$ be the convex cone generated by $P_\rho'$, that is, the set of sums of elements of $P_\rho'$.
Observe that if $\rho_1$ and $\rho_2$
are equivalent representations, then $P'_{\rho_1}=P'_{\rho_2}$ and $P_{\rho_1}=P_{\rho_2}$.

Let $G$ be a topological group and $\mathcal R$ a set of equivalence classes of unitary representations of $G$.
The \emph{Fell topology} on $\mathcal R$ is defined as follows: a typical neighborhood of $[\rho]\in \mathcal R$ has the form
$$
W([\rho];\{f_1, \cdots, f_n, C, \gep\})=\{[\s]\in\mathcal R :
\, \exists g_1, \cdots, g_n\in P_\s
\  \forall x\in C  \
|f_i(x)-g_i(x)| <\gep\},
$$
where $f_1, \cdots, f_n\in P_\rho$ (or $P_\rho'$), $C$ is a compact subspace of $G$, and $\gep>0$.
In particular, the Fell topology is defined on the dual object $\widehat G$.

Denote by $\repng$ the set of all
continuous $\mbox{$n$-dimensional}$ unitary representations of
a topological group $G$, i.~e., the set of all continuous homomorphisms of $G$ into
the unitary group $\mathbb U(n)$, equipped with the compact-open topology.
It follows from a result due to Goto \cite{goto} that $\repng$ is a
locally compact and uniformizable space (the space $\repg=\sqcup
_{n<\omega }\repng$ (as a topological sum) is called the
\textit{Chu dual} of $G$ \cite{chu}).

The set $G_n^x=\frac{\repng}{\sim}$ of equivalence classes of
unitary representations of dimension $n$ may be equipped with the quotient
topology induced by the mapping
$$Q_n:\repng \longrightarrow  G_n^x.$$
It is easily seen that on each $G_n^x$ this quotient topology
coincides with the Fell topology defined above. Furthermore, if
$G$ is a compact group then each partial dual space $G_n^x$ is uniformly
discrete.

We also equip the subset $\irn$, of all irreducible representations of
dimension $n$, with the compact open topology inherited from $\repng$ and
we denote by $\widehat{G}_n=\frac{\irn}{\sim}=Q_n(\irn)$ the
quotient space formed by the equivalence classes of irreducible
representations of dimension $n$. Again, the space $\widehat{G}_n$ is equipped with the Fell topology defined above.
Thus the partial dual object $\widehat{G}_n$ is a topological subspace of $ G_n^x$ for all $n\in\NN$.

If $G,H$ are groups, then $Hom(G,H)$ is the set of all homomorphisms from $G$ to $H$.
If $G,H$ are topological groups, then $CHom(G,H)$ is
the set of continuous elements of $Hom(G,H)$.

A (group) compactification of the topological group $G$ is a pair $(X,\phi )$, such that
$X$ is a compact group, $\phi\in CHom(G,X)$, and ${\rm ran}(\phi)$ is dense in $X$.

If $(X,\phi)$, $(Y,\psi)$ are two compactifications of $G$, then $(X,\phi)\leq (Y,\psi)$ if and only if
$\phi=\Gamma\circ \psi$\,  for some $\Gamma\in CHom(Y,X)$.  $(X,\phi)$ and $(Y,\psi)$ are \emph{equivalent}
if $(X,\phi)\leq (Y,\psi)$ and $(Y,\psi)\leq (X,\phi)$.

Observe that the $\Gamma$ in the definition above is unique and onto (since the range of
a compactification is dense). In the case that $(X,\phi)$ and $(Y,\psi)$ are equivalent, $\Gamma$ is a continuous isomorphism.

The \emph{Bohr compactification} of an  arbitrary topological group $G$
is defined as the group compactification $(bG,b)$
with the following universal property:
for every continuous homomorphism $h$
from $G$ into a compact group $K$ there is a continuous homomorphism $h^{b}$
from $bG$ into $K$ extending $h$
in the sense that $h=h^b \circ b$.
The Bohr compactification $(bG,b)$ is the unique (up
to equivalence) largest compactification of $G$ in the order $\leq$.

\section{Self-bohrifying groups}
If $G$ is a topological group, then $G_d$ denotes the same group
with the discrete topology. According to Hart and Kunen \cite{HartKunen:2001},
a compact group X is \emph{self-bohrifying} if
$(X,b)$ is the Bohr compactification of $X_d$, where $b\colon X_d\to X$ is
the identity map. That is $b(X_d)=X$. Previously, in [6, page 195], Comfort had defined a topological group $G$
to be a \emph{van der Waerden group} ($vdW$-group) if every homomorphism from $G$ to a compact group is continuous.
It is easily seen that a compact group is a $vdW$-group if and only if its topology is the finest
totally bounded group topology on $G$.
Finally, a topological group $G$ is \emph{tall} if for each positive integer $n$
there are only finitely many classes $[\rho]\in \widehat G_n$.
It is remarkable that
all the notions introduced above are related for compact groups (cf. \cite[Th. 1.15]{HartKunen:2001}).

\bprp\label{p:1}
Let $G$ be a compact group and let $G_d$ be the same algebraic group equipped with the discrete topology.
The following conditions are equivalent:
\begin{enumerate}
\item $G_d$ is tall.
\item $G$ is self-bohrifying.
\item $\widehat{(G_d)}_n\cup\{1_G\}$ is discrete for all $n\in\N$.
\item $G$ is a $vdW$-group.
\end{enumerate}
\eprp
\bpf
$(1)\Rightarrow (3)$ is obvious.

$(3)\Rightarrow (1)$. If $G_d$ is not tall, then some partial dual $\widehat{(G_d)}_n$ must be infinite.
Applying Proposition 5.5 of \cite{FHU}, it follows that $\widehat{(G_d)}_m\cup\{1_G\}$ is not discrete
for some $1\leq m\leq n^2$.

$(1)\Leftrightarrow (2)$ can be found in \cite{HartKunen:2001}.

$(2)\Leftrightarrow (4)$ is clear.
\epf

For profinite groups,  it may be more useful to look at the smaller subclass of  representations with finite range.
Indeed, suppose that $G$ is a profinite group with a
neighbourhood base at the identity $\{N_i: i\in I\}$ consisting of open normal subgroups
of $G$, then $G$ is the projective limit of the family $\{G/N_i : i\in I\}$ consisting of finite groups.
Each representation with finite range $\rho\in \widehat{G}^f$ has a kernel containing some $N_i$.
Therefore we have
$$\widehat{G}^f=\bigcup\limits_{i\in I} \widehat{\left(\frac{G}{N_i}\right)},$$
where we identify the representations of a quotient of $G$ with representations of $G$.
If $G$ is a residually finite group, we denote by $\overline G$ its profinite completion, which is a group
compactification of $G$.

Next result is a variant of Proposition \ref{p:1} for profinite groups.

\bprp
Let $G$ be a profinite group and let $G_d$ the same algebraic group with the discrete topology.
The following conditions are equivalent:
\begin{enumerate}
\item $G_d$ is tall for representations with finite range; that is, the set $(\widehat{G_d})^f_n$, of equi\-va\-lence classes of $n$-dimensional
irreducible representations with finite range of $G_d$ is finite for each positive integer $n$.
\item $\overline{G_d}=G$.
\item $(\widehat{G_d})^f_n\cup\{1_G\}$ is discrete for all $n\in\N$.
\end{enumerate}
\eprp
\bpf
$(1)\Longrightarrow (2)$. By definition  $\overline{G_d}=G$ if and only if every homomorphism
of finite rank on $G$ is continuous. Moreover, the latter is equivalent to that every abstract subgroup
of finite index in $G$ is open. On the other hand, $G_d$ is tall for representations with finite range if and only if it has only finitely many abstract
subgroups of index n for each integer n. Thus, by \cite[Th. 2]{smith_wilson:03}, if $G_d$ is tall for representations with finite range, every
abstract subgroup of finite index is open, which yields $(2)$.

$(2)\Longrightarrow (3)$. If $(3)$ is false, then for some $n\in \N$, the set $(\widehat{G_d})^f_n$ is infinite.
According to \cite[Th. 2]{smith_wilson:03}, this means that some abstract subgroup of finite index is not open,
which means that there are group homomorphisms of finite range defined on $G$ that are not continuous.
Thus $(2)$ is also false.

$(3)\Longrightarrow (1)$. If $G_d$ is not tall for representations with finite range, it follows that the set $(\widehat{G_d})^f_n$ is infinite for some
$n\in \N$. Applying Proposition 5.5 of \cite{FHU}, it follows that $(\widehat{G_d})^f_m\cup\{1_G\}$ is not discrete
for some $1\leq m\leq n^2$.
\epf
\bigskip

\section{Bohr compactification versus Chu duality}
In the previous secction we have seen how selfborifying groups are characterized by the property that their partial duals $\widehat{G}_n$
are finite for all $n\in \NN$. In this section, we explore when the Bohr compactification of a locally compact group $G$ coincides with its Chu or unitary bidual.
The main feature of Chu duality is the construction of a \emph{bidual} of $G$
from $\repg$. This bidual is formed by the so-called quasi-representations. If we define $\mathbb{U}=\sqcup _{n<\omega
}\mathbb{U}(n)$ (topological sum), a
\textit{quasi-representation} of $\repg$ is a continuous mapping
$p:\repg\longrightarrow \mathbb{U}$ which preserves the main operations defined
between unitary representations: direct sums, tensor products,
unitary equivalence and sends the elements of $\repng$ into
$\U(n)$ for all $n\in \N$ (see \cite{chu} or \cite{heyer70} for details).
The set of all quasi-representations of $\repg$ equipped
with the compact-open topology on $\repg$ is a topological group with
pointwise multiplication as the composition law, called the
\textit{Chu bidual group} of $G$ and denoted by $G^{xx}$.
Thus, a \nbd base of the identity of $G^{xx}$ consists of sets of
the form $[K_n,V]=\{p\in G^{xx} : p(K_n)\subset V \}$, where $V$
is any \nbd of the identity in $\U(n)$ and $K_n$ is any compact
subset of $\repng$, $n\in \N$. It is easily verified that the
evaluation map $\mathcal E_{G} :G\longrightarrow G^{xx}$ is a group
homomorphism which is a monomorphism if and only if $G$ is MAP.
We say that $G$ \textit{satisfies Chu duality}
when the evaluation map $\mathcal E_{G}$ is an isomorphism of
topological groups. In this terminology, it was shown in \cite{chu}, that LCA groups
and compact groups satisfy Chu duality (indeed Chu duality reduces
to Pontryagin duality and to Tannaka duality, respectively, for such
groups). Here, the group $\mathcal E_{G}(G)$ is always assumed to be
equipped with the topology inherited from $G^{xx}$.
Remark that $G^{xx}_{|\repng}\subseteq C(\repng, \mathbb{U}(n))$ ,
where $C(X,\mathbb U(n))$ denotes the space of all continuous
functions from a topological space $X$ into $\mathbb U(n)$.

Our first result  demonstrates the existence of non compact Chu reflexive groups whose partial dual spaces $\widehat G_n$,
equipped with the Fell topology, are discrete.
Theorems \ref{t:01} and \ref{t:1} fix a wrong assertion stated in \cite[Theorem 4.4]{galihern04} (see also \cite[Proposition 2.2]{her_wu_2006}).
\bigskip

\bthm\label{t:01}
Let $H$ be a finite simple non abelian group, let $I$ be an infinite index set and consider the discrete group $G=H^{(I)}$.
Then $G$ is Chu reflexive and $(\widehat{G}_n,t_p(G))$ is a discrete topological space for all $n\in\NN$.
\ethm

\bpf Remark that the group $G$ also accepts a precompact  topology as a dense subgroup of the compact group $H^I$.
We want to prove that $(\widehat{G}_n,t_p(G))$ is a discrete topological space.
To see this, we show that for every $\varphi\in\irn$ there is a finite subset $F$ of $G$ and $\epsilon>0$ such that if
$\psi\in\irn$ and $\psi\in N(\varphi,F,\epsilon)$, then $\psi\sim\varphi$.

This is vacuously true for $n=1$ because $H$ is a finite simple non abelian group.

For $n>1$, note that if $\varphi\in \textrm{irrep}_n(G)$, then $\varphi\sim\varphi_{i_1}\otimes\cdots\otimes\varphi_{i_N}$,
$\varphi_{i_l}\in \textrm{irrep}_{n_l}(H)$, $n_l>1$, $1\leq l\leq N$ and $n_1\cdots n_N=n$.
Set $F_{i_1,\cdots,i_N}=\prod\limits_{i\in I}B_i$, where $B_i=\{e\}$, $i\neq i_l$, and $B_{i_l}=H$, $1\leq l\leq N$.
We proceed by induction over $N$.

Suppose first that $\varphi(x)=(\varphi_{i_1}\circ \pi_{i_1})(x)$, $x=(x_i)_{i\in I}\in G$,
where $\pi_{i_1}\colon G\to H_{i_1}$ is the canonical projection of $G$ onto $H_{i_1}=H$.
With some notational abuse, we identify $\varphi|_{F_{i_1}}$ with $\varphi_{i_1}$. Therefore we have
$\varphi|_{F_{i_1}}=\varphi_{i_1}\in \textrm{irrep}_n(H)$.
Take $F=F_{i_1}$, which is a finite subgroup of $G$, and $\epsilon=\frac{1}{n}$.
It is known that if $\psi\in\textrm{irrep}_n(G)$ and $\|\varphi(x)-\psi(x)\|<\frac{1}{n}$ for all $x\in F=F_{i_1}\cong H$,
then $\varphi_{i_1}=\varphi|_{F_{i_1}}\sim\psi|_{F_{i_1}}$ (see \cite[Satz 12.1.6]{heyer70}).
Moreover, assuming that $\psi\sim\psi_{j_1}\otimes\cdots\otimes\psi_{j_M}$, $\psi_{i_l}\in \textrm{irrep}_{m_l}(H)$, $1\leq l\leq M$
and $m_1\cdots m_M=n$. If $i_1\notin\{j_1,\cdots,j_M\}$ then $\psi|_{F_{i_1}}\sim I_n$, where $I_n$ is the $n$-dimensional identity matrix.
Therefore $\varphi_{i_1}\sim I_n$, which is impossible. Thus, there is an index $j_k=i_1$.
Hence  $\varphi|_{F_{i_1}}\sim\psi|_{F_{i_1}}$ yields
$$\varphi_{i_1}=\varphi|{F_{i_1}}\sim\psi|{F_{i_1}}\sim\psi_{i_1}\otimes I_{\frac{n}{m_k}}\sim
\psi_{i_1}\oplus\underbrace{\cdots}_{\frac{n}{m_k}}\oplus \psi_{i_1}.$$
Since, the decomposition as a direct sum of irreducible representations is unique except for order and equivalence (see \cite{heyer70}),
we deduce $\frac{n}{m_k}=1$ and $\varphi_{i_1}\sim\psi_{i_1}$. That is
$$\psi=\psi_{i_1}\circ \pi_{i_1}\sim\varphi_{i_1}\circ \pi_{i_1}=\varphi.$$
Thus, the case $N=1$ is verified.

Assume now that the assertion holds for $N$ and consider
$$\varphi=(\varphi_{i_1}\otimes\cdots\otimes\varphi_{i_N}\otimes\varphi_{i_{N+1}})\circ\pi_{(i_1,\dots ,i_{N+1})},$$ where
$$\pi_{(i_1,\dots ,i_{N+1})}\colon G\to H_{i_1}\times\dots H_{i_{N+1}}$$ is the canonical projection,
$\varphi_{i_l}\in \irs_{n_l}(H)$, $1\leq l\leq N+1$ and $n_1\cdots n_{N+1}=n$.

We now look at $\varphi$ as follows: $$\varphi = (\varphi_{1}\otimes\varphi_{2})\circ\pi_{(1,2)},$$ where
$$\varphi_1 =\varphi_{i_1}\otimes\cdots\otimes\varphi_{i_N},$$
$$\varphi_2=\varphi_{i_{N+1}},$$
$$\pi_{(1,2)}\colon G\to H^{(I\setminus\{i_{N+1}\})}\times H_{i_{N+1}}$$ is the canonical projection,
$\varphi_1\in\irs_{m_1}(H^{(I\setminus\{i_{N+1}\})})$, $m_1=n_1\dots n_N$, and $\varphi_{i_{2}}\in\irs_{n_{N+1}}(H)$.

Again, with some notational abuse we have
$$\varphi|_{H^{(I\setminus\{i_{N+1}\})}}=\varphi_{1}\in\textrm{irrep}_{m_1}(H^{(I\setminus\{n_{N+1}\})})$$ and
$$\varphi|_{F_{i_{N+1}}}=\varphi_2\in \textrm{irrep}_{n_{N+1}}(H).$$
Take $F=F_{i_1,\cdots,i_{N+1}}$, which is a finite subgroup of $G$, and $\epsilon=\frac{1}{n}$.
Let $\psi\in\textrm{irrep}_n(G)$ such that $\|\varphi(x)-\psi(x)\|<\frac{1}{n}$ for all $x\in F$.
 Then $$\|\varphi(x)-\psi(x)\|<\frac{1}{n}\ \text{for all}\ x\in F_{i_1,\cdots,i_{N}}$$ and
 $$\|\varphi(x)-\psi(x)\|<\frac{1}{n}\ \text{for all}\ x\in F_{i_{N+1}}.$$

 By our induction hypothesis, and using \cite[Satz 12.1.6]{heyer70}, we deduce that
$$\varphi_1=\varphi|_{H^{(I\setminus\{i_{N+1}\})}}\sim\psi|_{H^{(I\setminus\{i_{N+1}\})}}=:\psi_1\in\textrm{irrep}_{m_1}(H^{(I\setminus\{i_{N+1}\})})$$
and
$$\varphi_{2}=\varphi|_{F_{i_{N+1}}}\sim\psi|_{F_{i_{N+1}}}=:\psi_2\in\textrm{irrep}_{n_{N+1}}(H).$$

Therefore
$$\psi_1\sim\varphi_{i_1}\otimes\cdots\otimes\varphi_{i_N}$$ and
$$\psi_2\sim\varphi_{i_{N+1}}.$$
This implies that $$\psi\sim \psi_1\otimes\psi_2\sim \varphi_{i_1}\otimes\cdots\otimes\varphi_{i_N}\otimes\varphi_{i_{N+1}} = \varphi.$$

It was proved in \cite{chu} that the discrete group  $G=H^{(I)}$ is Chu reflexive. Nevertheless, we include a short version
of it here for the reader's sake. Let $\phi\in\repnh$ be a faithful unitary representation of $H$ and set
$\phi_i=\phi\circ\pi_i$, where $\pi_i\colon G\to H_i$ is the canonical projection of $G$ onto $H_i=H$.
Since $G=H^{(I)}$ is a discrete group, it is easyly verified that the trivial representation $I_n$
(the one that sends every element to $I_n$) is a closure point of the discrete subset
$\{\phi_i : i\in J\}$ for every infinite subset $J\subseteq I$.
On the other hand, we have that $bG\cong H^I$ (see \cite[Cor. 3.10]{chu}). Thus, ir order to show that $G$ is Chu reflexive,
it will suffice to verify that no element in $H^I\setminus G$ defines a continuous quasi-representation.
Hence, take $(x_i)\in H^I\setminus G$. Then there is an infinite subset $J\subseteq I$ such that
$x_i\not= e_H$ for all $i\in J$. Furthermore, since the group $H$ is finite, we may assume without loss of generality
that $x_i=x_j$ for all $i,j$ in $J$. As a consequence, there is $A\in \mathbb{U}_n$ such that
$(x_i)[\phi_i]=\phi(x_i)=A\not=I_n$ for all $i\in J$. This implies that $(x_i)$ is not continuous on $\repng$,
which completes the proof.
\epf
\bigskip

The following result completes Theorem \ref{t:01} by giving a characterization of locally compact groups whose Chu quasi-dual coincides
with  their Bohr compactification.

\bthm\label{t:1}
Let $G$  be a locally compact group. It holds that $G^{xx}$ is
topologically isomorphic to $bG$ if and only if the space $G_n^x$ is
discrete for any $n\in \N$.
\ethm
\bpf
{\it Necessity:} Let $\mathcal E_n : \repng\longrightarrow {\rm rep}_n(G^{xx})$ the canonical evaluation map defined by
$\mathcal E_n(D)(p)=p(D)$ for all $D\in \repng$ and $p\in G^{xx}$.
We first check that $\mathcal E_n$ is continuous and injective. For the continuity,
taking into account that the spaces $\repng$ and ${\rm rep}_n(G^{xx})$ are equipped with the compact open topology
on $G$ and $G^{xx}$, respectively, let $\{D_j \}$ be a net converging to $D$ in $\repng$.
Given an arbitrary compact subset $K$ of $G^{xx}$, since $\repng$ is locally
compact and $G^{xx}_{|\repng}\subseteq C(\repng,\mathbb U(n))$,
it follows that $K$ is an equicontinuous subset of
$C(\repng,\mathbb U(n))$. Furthermore, the net $\{D_j \}$ converges to $D$ in the compact open topology on $\repng$
and every element $p\in G^{xx}$ is continuous on $\repng$, which means that the net $\{p(D_j) \}$ converges to $p(D)$
for all $p\in G^{xx}$. That is to say, the net $\{\mathcal E_n(D_j) \}$ converges to $\mathcal E_n(D)$ pointwise on $G^{xx}$.
Hence, since $K$ is an equicontinuous subset of $C(\repng,\mathbb U(n))$, it follows that $\{\mathcal E_n(D_j) \}$ converges to
$\mathcal E_n(D)$ uniformly on $K$. This yields the continuity of $\mathcal E_n$.

For the injectivity, suppose that $D,E$ belong to
$\repng$ and $\mathcal E_n(D)=\mathcal E_n(E)$. Then $p(D)=p(E)$ for all
$p\in G^{xx}$. In particular, $D(g)=E(g)$ for all $g\in G$. Thus, D=E.

Now, we define $\Delta_n : G^x_n\longrightarrow (G^{xx})_n^x $ \ so that the
following diagram commutes
\[
\begin{CD}
\repng@>\mathcal E_n>>   {\rm rep}_n(G^{xx})\\
@VQ_nVV   @VVQ_n^{xx}V\\     
G^x_n @>\Delta_n>>  (G^{xx})_n^x
\end{CD}
\]
\medskip

\noindent where $Q_n$ and $Q_n^{xx}$ are the canonical quotient
mappings.

Since we are assuming that $G^{xx}$ is topologically isomorphic to $bG$,
it follows that $G$ is dense in $G^{xx}$ and, using this fact,
it is easy to see that $\Delta_n$ is well defined and $1$-to-$1$. The fact that $Q_n$ and
$Q_n^{xx}$ are quotient mappings yields the continuity of
$\Delta_n$.
On the other hand, since $G^{xx}$ is compact,
it follows  that $(G^{xx})^x_n$ is discrete (see \cite{heyer70}).
Since $\Delta_n$ is continuous and injective, it follows that
$G^x_n$ is also discrete. This completes the proof.

{\it Sufficiency:}
The Bohr compactification $bG$ of $G$ can be realized as the set of {\it all}
(con\-ti\-nuous or not) quasi-representations of $\repg$ endowed with the topology of pointwise
convergence. On the other hand, the group $G^{xx}$ is complete
and there is a con\-ti\-nuous isomorphism of $G^{xx}$ onto a dense subgroup of $bG$
(see \cite{galihern04}). Hence, in order to prove that $bG\cong G^{xx}$, it suffices to prove that $G^{xx}$ is
precompact.\
To do this, we will verify that, on $G^{xx}$, the compact open topology and the pointwise topology on $\repg$
coincide. Remark that the
latter topology is trivially precompact since embeds $G^{xx}$ into a product of
finite dimensional unitary groups.

Therefore, let $K$ be an arbitrary compact subset of $\repng$. The continuity of $Q_n$ implies that
$Q_n(K)$  will be a compact in $G^x_n$ and, by hypothesis, necessarily  finite.
Form a finite subset  $\Phi_n=\{\sigma_1,\ldots ,\sigma_k\}\subseteq \repng$,
taking a representative of each element in $K$.
The unitary groups $\mathbb U(n)$ contain a neighborhood basis of the identity $I_n$ whose elements are invariant
under conjugation. Furthermore, the density of $G$ in $G^{xx}$ implies that the elements in the latter group preserve conjugation.
Altogether, this yields that the convergence on the finite set $\Phi_n$ is the same as uniform convergence on $K$, which completes the proof.
\epf
\bigskip

The next result improves Theorem \ref{t:1} for discrete groups, clarifies Theorem \ref{t:01}, and extends Proposition \ref{p:1}.

\bcor\label{t:2}
Let $G$  be a discrete group. Then $G^{xx}$ is
topologically isomorphic to $bG$ if and only if the space $\widehat{G}_n$ is
finite for any $n\in \N$.
\ecor
\bpf
{\it Necessity:} Since $G$ is discrete, the dual object $\widehat{G}$ is quasicompact (see \cite[$\S$18.1]{dixmier}).
On the other hand, the partial dual $\widehat{G}_n$ is a topological subspace of $G^x_n$ that is discrete by Theorem \ref{t:1}.
As $\bigcup\limits_{i=1}^{n}\widehat{G}_i$ is closed in $\widehat{G}$, it is also quasicompact and discrete, hence finite.

{\it Sufficiency:} As in Theorem \ref{t:1}, in order to prove that $bG=G^{xx}$ it suffices to prove that $G^{xx}$ is
precompact. Now, since $\widehat{G}_m$ is finite, it follows that the quotient uniformity canonically associated to $\widehat{G}_m$
by the quotient mapping $Q_{m}$ is the discrete uniformity. Hence,
the subset  $\Lambda_n :=\{\sigma_i\}_{i\in I_n} \subset \irn$,
formed by choosing a representative of each equivalence class, is finite and, hence, inherites the discrete uniformity
from $\rep_n(G)$. We next show how the topology of uniform convergence on compact subsets
of $\rep_n(G)$ and the topology of pointwise convergence on $\Lambda_n$ coalesce. Then the
proof is completed by observing that the latter topology is trivially precompact.

Indeed, choose  an arbitrary compact subset  $K\subseteq \repng$.
Each element of $K$ decomposes as the direct sum of a finite number
of irreducible subrepresentations, let's denote by $[K]_m\subseteq
\irm$ the subset of
all irreducible components of elements of $K$
of dimension $m$. By our initial assumption, we have that
$\Q_m([K]_m)$ is a finite subset of $\widehat{G}_m$ for $1\leq m\leq n$.
Form a finite subset $\{\sigma_1,\ldots ,\sigma_k\}$ of $\Lambda$
taking a representative of each element in $\Phi_n=\bigcup_{m=1}^{n}\Psi_m([K]_m)$.
Arguing as in the last paragraph of Theorem \ref{t:1},
we conclude that pointwise  convergence on the finite set $\Phi_n$ is the same as uniform convergence on $K$.
\epf
\bigskip

\brem\label{rem1}
It was stated in \cite[Theorem 4.4]{galihern04} that if $G$ is a locally
compact group such that $\widehat G_n$ is discrete for all $n\in \N$,
then $bG=G^{xx}$ and therefore the group $G$ may not be Chu reflexive.
This result is incorrect as we have shown in Theorem \ref{t:01} and
the reason why this happens is subtle. It was proven in \cite[Cor. 5.6]{FHU}
that if $G$ is a discrete group and some partial dual object $\widehat G_n$ is infinite, then the trivial representation
$1_G$ is not an isolated point in in some $\widehat G_m\cup\{1_G\}$ with $1\leq m\leq n^2$.
This means that although the topological space $\widehat G_m$ may be discrete, the \emph{uniformity} that $\widehat G_m$
inherits from $\widehat G$ is not, as asserted in \cite{galihern04}.
This said, we want to emphasize that the main consequences of \cite[Theorem 4.4]{galihern04}
are preserved without mo\-di\-fi\-cation using Theorem \ref{t:1}.
In particular, we have:\medskip
\erem

\bcor
If $G$ is an  infinite vdW-group, equipped  with
the discrete topology, then $G$ does
not satisfy Chu-duality.
\ecor
\bcor
A Kazhdan group satisfies Chu duality if and only if it is compact.
\ecor

\section{Acknowledgement}

\noindent {We thank the referee for her/his helpful comments.}

\end{document}